\documentclass[a4paper,12pt]{article}
\usepackage{amsfonts}
\usepackage{graphicx}
\usepackage{amssymb}
\usepackage[top=1.6in, bottom=1.5in, left=1.16in, right=1.16in]{geometry}
\newtheorem{theo}{Theorem}[section]

\newtheorem{lem}{Lemma}[section]

\newtheorem{defi}{Definition}[section]
\newcommand{\be}{\begin{equation}}
\newcommand{\ee}{\end{equation}}
\newcommand\bes{\begin{eqnarray}}
\newcommand\ees{\end{eqnarray}}
\newcommand{\bess}{\begin{eqnarray*}}
\newcommand{\eess}{\end{eqnarray*}}

\begin{document}
\setlength{\baselineskip}{1.1\baselineskip} \pagestyle{myheadings}

\title{ \bf  \Large Dynamical behavior of a harvest single species model on growing habitat\thanks{The work is partially supported by PRC grant NSFC
(61103018) and NSF
 of Jiangsu Province (BK2012682,10KJB110011).}}
\date{\empty}
\author{\normalsize Zhi Ling$^{a}$,  Lai Zhang$^{b}$\thanks{Corresponding author: zhling@yzu.edu.cn}  \\
{ \small $^{a}$ School of Mathematical Science, Yangzhou University,}\\
{ \small Yangzhou 225002, P. R. China}\\
{ \small $^{b}$ Department of Mathematics and Mathematical
Statistics, Ume{\aa} University,}\\
{ \small SE-90187, Ume{\aa} Sweden}} \maketitle

\begin{quote}
\noindent {\bf Abstract:} {This paper is concerned with a
reaction-diffusion single species model with harvesting on
$n$-dimensional isotropically growing domain. The model on growing
domain is derived and the corresponding comparison principle is
proved. The asymptotic behavior of the solution to the problem is
obtained by using the method of upper and lower solutions. The results
show that the growth of domain takes a positive effect on the
asymptotic stability of positive steady state solution while it
takes a negative effect on the asymptotic stability of the trivial
solution, but the effect of the harvesting rate is opposite. The
analytical findings are validated with the numerical simulations.}


\end{quote}

\section{Introduction}\label{s:1}
The growth of a single species population that evolves according to
a logistic law while assuming that species undergoes a random walk
may be modeled by the following equation
\begin{eqnarray}
\label{10}\left.
\begin{array}{lll}
\displaystyle u_t=d\Delta u+r u(1-\frac{u}{K}),&x\in\Omega,t>0,
\end{array} \right.
\end{eqnarray}
where $d$ denotes the diffusion coefficient, $r$ represents the
intrinsic growth rate and $K$ is the natural carrying capacity of
the environment. We take these parameters to be positive constants.
$\Omega$ is a bounded subset of $\mathbb{R}^n$ ($n\geqslant1$) with
smooth boundary $\partial\Omega$. And the environment $\Omega$ is
homogeneous (i.e., the diffusion does not depend on $x$). $u(x,t)$
is the density of the species at position $x$ and time $t$. $u_t=\partial u/\partial t$, $\Delta$ denotes the Laplace operator in $\Omega$. Equation
(\ref{10}) is often called Fisher's equation after Fisher \cite{F},
who proposed the one-dimensional version as a model for the spread
of an advantageous gene in a population, and it was also studied by
Kolmogoroff, Petrovsky and Piscounoff \cite{KPP}, who studied the
equation in depth and obtained some of the basic analytical results.

We assume that the species migrates in a domain surrounded by a
hostile environment, so we can consider the initial and boundary
conditions as
\begin{eqnarray}
\label{10'}\left.
\begin{array}{lll}
u(x,t)=0,&x\in\partial\Omega,\ t>0,\\
u(x,0)=u_0(x)\geqslant0 ,&x\in\Omega,
\end{array} \right.
\end{eqnarray}
where $u_0\in C^2(\bar{\Omega})$ and $u_0=0$ on $\partial\Omega$.
The dynamics of (\ref{10}) and (\ref{10'}) has been completely
studied, see \cite{DH}. Here we briefly describe the results which
are related to present paper. Let $\lambda_1$ be the principal
eigenvalue of the problem
\begin{eqnarray*}
\label{101}\left\{
\begin{array}{lll}
-\Delta \phi=\lambda\phi,& x\in\Omega,\\
\phi(x)=0 ,& x\in\partial\Omega,
\end{array} \right.
\end{eqnarray*}
then we state the following Theorem.
\begin{theo}\label{T11}
For (\ref{10}) and (\ref{10'}), the following facts hold:\\
(1) If $0<r\leqslant d\lambda_1$, then there is only one nonnegative
steady state solution $u=0$,
which is globally asymptotically stable, that is, for any nonnegative nontrivial $u_0$, $\lim_{t\to \infty} u(x,t)=0$ uniformly.\\
(2) If $r>d\lambda_1$, then there is only one positive steady state
solution $u=u^*(x)$, which is globally asymptotically stable, that
is, for any nonnegative nontrivial $u_0$, $\lim_{t\to \infty}
u(x,t)=u^*(x)$ uniformly.
\end{theo}

 From the point of view of human needs, the exploitation of
biological resources and the harvest of population are commonly
practiced in fishery, forestry and wildlife management. Concerning
the conservation for the long-term benefits of humanity, there is a
wide-range of interest in the use of bioeconomic modeling to gain
insight in the scientific management of renewable resources like
fisheries and forestry. At the same time, harvesting has a strong
impact on the dynamic evolution of a population. Hence it is natural
to add the harvesting term to the right-hand side of the first
equation in (\ref{10}), and the equation would be
\begin{eqnarray}
\label{110}\left.
\begin{array}{lll}
\displaystyle u_t=d\Delta u+r u(1-\frac{u}{K})-hu,
\end{array} \right.
\end{eqnarray}
 where $h>0$ is a parameter which represents
the level of harvesting, $hv$ is the harvesting yield per unit time.
One can see Murray \cite{M} for details about the harvest model.

As we know that the conventional theory of harvested populations
basing on equations in which the various environmental is treated as
fixed domain. In fact, the ecological environment is not always the
same in nature, the habitats of species usually changes due to many
reasons, for example, some insects live on a growing leaf, some
fishes live in an expanding river due to a warming effect, some
animals live in desert which is expanding continuously. A natural
question arises that how species react to the changing of their
habitats. Take this into account, in present paper, we will consider
the problem on growing domain.

Indeed, domain growth has been suggested as an important mechanism
in pattern formation and election, we refer to \cite{CGM,CHM,HMS,
SG, MM,11,JM,CAB} and the references therein for more details.
However, since the presence of time-dependent transport coefficients
in the equations which constructed on growing domain leads to
difficulty in stability analysis, most of known work was carried out
though numerical computation and simulations.

Recently, Tang etc \cite{TL} considered a diffusive logistic
equation on one dimensional isotropically growing domain with linear
growth function and exponential growth function respectively and get
the asymptotic behavior of the solution by constructing upper and
lower solutions. In this paper, we try to use this method to study
the asymptotic behavior of solution to problem (\ref{10'}) and
(\ref{110}) on $n$ dimensional growing domain.

The organization of this paper is as follows: In Section 2, a general
reaction-diffusion equation with domain growth is developed in
$n$-dimensional space $\mathbb{R}^n$ and then the harvest single
species logistic model on an isotropically growing domain is
constructed. In Section 3, we restrict our attention to the
isotropically growing domain and analyse the asymptotic behavior of
solutions. In section 4, by performing a series of simulations, we
illustrate our analytical result. Finally, we give a brief
conclusion in Section 5.

\section{Model on growing domain}\label{s2}\setcounter{equation}{0}
In this section, we first model a general reaction-diffusion
equation on growing domain in $\mathbb{R}^n$ and then present the
single species harvest model on an isotropically growing domain. The
approach is as in \cite{CGM}.

 Let $\Omega (t)\subset \mathbb{R}^n$
be a simply connected bounded growing domain at time $t\geqslant 0$
with its growing boundary $\partial\Omega(t)$.  For any point
$x(t)=(x_1(t),x_2(t), \ldots, x_n(t)) $ $ \in \Omega(t)$, we assume
that $u(x(t),t)$ is the density of a species, at position $x(t)$ and
time $t\geqslant 0$. According to the principle of mass
conservation, we have
 \bess
    \frac{\mathrm{d}}{\mathrm{d}t}\int_{\Omega(t)} u(x(t),t) \mathrm{d}x
    = -\int_{\partial\Omega(t)}\textbf{\emph{J}} \cdot \textbf{\emph{n}} \mathrm{d}S
    +\int_{\Omega(t)} f(u)\mathrm{d}x,
\label{b01}
 \eess
where $\textbf{\emph{J}}$ is the flux across the boundary
$\partial\Omega(t)$, $\textbf{\emph{n}}$ is the outward normal
vector on $\partial\Omega(t)$, $f(u)$ is the reaction term within
the domain. Using the divergence theorem, the above equation becomes
 \bes
    \frac{\mathrm{d}}{\mathrm{d}t}\int_{\Omega(t)} u(x(t),t) \mathrm{d}x
    =\int_{\Omega(t)} [-\nabla \cdot {\textbf{\emph{J}}}
     +f(u)]\mathrm{d}x.
 \label{b02} \ees

On the other hand, the growth of domain generates a flow velocity
field $\textbf{\emph{a}}=(\dot x_1(t),\dot x_2(t),\ldots, \dot
x_n(t))$. Using the Reynold transport theorem to the left-hand side
of (\ref{b02}) yields
 \bess
    \frac{\mathrm{d}}{\mathrm{d}t}\int_{\Omega(t)} u(x(t),t)\mathrm{d}x
    =\int_{\Omega(t)}\left[\frac{\mathrm{d}u}{\mathrm{d}t} +u(\nabla\cdot \textbf{\emph{a}})\right]\mathrm{d}x,
 \label{b03}\eess
where $\frac{\mathrm{d}u}{\mathrm{d}t}$ is the total derivative of $u$, i.e. \bess
    \frac{\mathrm{d}u}{\mathrm{d}t}=\frac {\partial u}{\partial t}+\nabla u\cdot \textbf{\emph{a}}.
\label{b04}\eess Hence we can write (\ref{b02}) as follow:
 \bess
    \int_{\Omega(t)} \left[\frac{\partial u}{\partial t}+\nabla u\cdot \textbf{\emph{a}}
     +u(\nabla\cdot \textbf{\emph{a}})\right]\mathrm{d}x
     =\int_{\Omega(t)}[-\nabla \cdot {\textbf{\emph{J}}}+f(u)]\mathrm{d}x.
  \label{b05}
\eess
 Since $\Omega(t)$ is arbitrary, then the differential equation
\bes
 \frac{\partial u}{\partial t}+\nabla u\cdot \textbf{\emph{a}}
 +u(\nabla\cdot \textbf{\emph{a}})= -\nabla \cdot {\bf \textbf{\emph{J}}}+ f(u) & \  \textrm{in}\  \Omega(t)
 \label{b06}\ees
holds for any $(x,t)$. Assume the species undergoes a random walk,
the diffusion flux of $u$ follows Fick's law:
 \bess{\textbf{\emph{J}}}=-d\nabla u,
 \label{b07}
 \eess
 where $d$ is the diffusive coefficient of $u$. Thus the equation
(\ref{b06}) becomes \bes
 \frac{\partial u}{\partial t}+\nabla u\cdot\textbf{\emph{a}}
 +u(\nabla\cdot \textbf{\emph{a}})= d\nabla ^2 u+ f(u) &  \  \textrm{in}\  \Omega(t),
 \label{b08}\ees
where $\nabla u \cdot \textbf{\emph{a}}$ is called advection term
while $(\nabla \cdot\textbf{\emph{a}}) u$ is called dilution term.

In most cases, it is difficult to study the properties of solution
to (\ref{b08}) because of the advection and dilution terms. Let
$y_1, y_2, \ldots, y_n$ be fixed cartesian coordinates in fixed
domain $\Omega(0)$ such that $x_1(t)=\hat x_1(y_1, y_2, \ldots, y_n,
t), x_2(t)=\hat x_2(y_1, y_2, \ldots, y_n, t), \ \dots, \
x_n(t)=\hat x_n(y_1, y_2, \ldots, y_n, t).$ As $t$ varies, the
coordinates $x_1, x_2, \ldots, x_n$ change position with time. These
positions are then mapped or transformed to a fixed position given
by the $y_1,y_2,\dots, y_n$ coordinates. Under this transformation,
we suppose $u$ is mapped into the new function defined as \bes
  u(x_1(t), x_2(t), \ldots, x_n(t), t)=v(y_1,y_2,\ldots, y_n, t).
 \label{b09}
  \ees
Thus the equation (\ref{b08}) can be translated to another form
which is defined on the fixed domain $\Omega(0)$ with respect to
$y=(y_1,y_2,\ldots, y_n)$. However, the new equation is also more
complicated \cite{mad2}. To further simplify  the model equations
(\ref{b08}), we assume that domain growth is uniform and isotropic,
that is, the growth of the domain takes place at the same proportion
in all directions as time elapses. In mathematical terms,
$x(t)=(x_1(t),x_2(t),\ldots, x_n(t))$ can be described as follow:
\bes
 x(t)=\rho(t)y, &  y\in \Omega(0),
\label{b10} \ees where $\rho(t)$ is called growth function subject
to $\rho(0)=1$ and $\dot \rho(t)\geqslant0$ for all $t>0$.

By (\ref{b09}) and (\ref{b10}), we have
$$\textbf{\emph{a}}=\dot x(t)
            =\dot\rho(t)y=\frac{\dot\rho}{\rho}x,$$ $$ v_t=u_t+\nabla u \cdot \textbf{\emph{a}},\qquad  \nabla \cdot\textbf{\emph{a}}=\frac{n \dot\rho}{\rho},\qquad \Delta
u=\frac{1}{\rho^2(t)}\Delta v,$$
 where $n$ is the dimension of space.
Then (\ref{b08}) becomes the following form
 \bess
 v_t=\frac{d}{\rho^2(t)}\Delta v
   -\frac{n\dot \rho(t)}{\rho(t)} v + f(v), & y \in \Omega(0),\,\ t>0.
\label{b11} \eess Then we obtain the following single species
harvest problem on the growing domain $\Omega(t)$: \bes\left\{
\begin{array}{ll}
  \displaystyle v_t=\frac{d}{\rho^2(t)}\Delta v
   -\frac{n\dot \rho(t)}{\rho(t)} v +rv(1-\frac{v}{K})-hv, & \ y \in \Omega(0),\ t>0,
\\
v(y,t)=0, & \  y\in \partial \Omega(0), \ t>0,\\
v(y,0)=u_0(x(0)),  & \ y\in \Omega(0).
\end{array}\right.
\label{b12} \ees

\section{Analysis of the asymptotic behavior} \label{s3}\setcounter{equation}{0}

In this section we will study the asymptotical behavior of the
solution of (\ref{b12}). Though there are many different kinds of
typical growth functions, such as linear growth, exponential growth,
logistic (or saturated) growth, in a phenomenological sense, the
logistic growth is a biologically reasonable growth function, see
\cite{pl} for more details. For this reason, we consider the following logistic growth function
\[
\rho(t)=\frac{\exp(kt)}{1+\frac{1}{m}(\exp(kt)-1)},
\]
where $k>0$ and $m>1$. Notice that $\rho(t)$ is continuously
differentiable on $[0,+\infty)$ and satisfies
\[\rho (0)=1,\ \ \dot\rho(t)>0,\ \   \lim_{t\to\infty}\rho(t)= m>1.
\]

Next we give the following definition of upper and lower solutions
of (\ref{b12}):
\begin{defi}\label{def1}
  A function
  $\tilde{v} \in C^{2,1}(\Omega(0) \times (0,\infty)) \cap C(\bar\Omega(0)\times [0,+\infty))$
 is called an upper solution of (\ref{b12}) if it satisfies
\bes\left\{
\begin{array}{ll}
 \displaystyle \tilde v_t \geqslant \frac{d}{\rho^2(t)} \Delta \tilde v-\frac{n\dot\rho(t)}{\rho(t)}\tilde v
 + r \tilde v(1-\frac{\tilde{v}}{K})-h\tilde{v},  & y \in \Omega(0), \, \ t>0,
\\
\tilde v(y,t)\geqslant 0,  & y \in \partial\Omega(0),\ t>0,\\
\tilde v(y,0)\geqslant v_0(y),  & y \in \Omega(0).
\end{array}\right.
\label{c01}\ees Similarly,  $\hat v(y,t) \in C^{2,1}(\Omega(0)
\times (0,+\infty)) \cap C(\bar\Omega(0) \times [0, +\infty))$ is
called a lower solution of (\ref{b12}) if it satisfies all the
reversed inequalities in (\ref{c01}).
 \end{defi}

To prove our main results, we recall the following two lemmas.
\begin{lem}\label{lem1}(Comparison Principle) Let $v(y,t)$ be a solution of (\ref{b12}) ,
 $\tilde v(y,t)$ and $\hat v(y,t)$ are upper and lower
solutions of (\ref{b12}) respectively, then $\hat v(y,t)\leqslant
v(y,t)\leqslant \tilde v(y,t)$ in $\bar \Omega(0) \times [0,
+\infty)$.
\end{lem}
{\bf Proof.} Define $w=\tilde{v}-v$, and it is easy to see that
$w(y,t)$ satisfies \bess\left\{
\begin{array}{ll}
\displaystyle w_t\geqslant \frac{d}{\rho^2(t)}\Delta
w-\frac{n\dot\rho(t)}{\rho(t)}w+ rw(1-\frac{\tilde{v}+v}{K})-hw,
 &y \in\Omega(0)  ,\ t>0,
\\
w(y,t)\geqslant0, & y\in\partial\Omega(0), \ t>0,\\
w(y,0)\geqslant 0,  & y\in\Omega(0).
\end{array}\right.
\label{l11} \eess Applying the maximum principle leads to
$$w(y,t)\geqslant0,\ \ y\in\Omega(0), \ t\geqslant0,$$
that is $\tilde{v}(y,t)\geqslant v(y,t)$, $ y\in\Omega(0),
t\geqslant0$. Similarly, $\hat{v}(y,t)\leqslant v(y,t)$ can be
proved.

\begin{lem} \label{lem2} Let $v(y, t)$ be a  nonnegative nontrivial solution of the
following problem
 \bess\left\{
\begin{array}{ll}
\displaystyle v_t= \frac{d}{\rho^2(t)}\Delta
v-\frac{n\dot\rho(t)}{\rho(t)}v+ rv(1-\frac{v}{K})-hv,
 &y \in\Omega(0)  ,\ t>0,
\\
v(y,t)=0, & y\in\partial\Omega(0), \ t>0,\\
v(y,0)=v_0(y)\geqslant 0,  & y\in\Omega(0).
\end{array}\right.
\label{c03} \eess If $v(y,0)\in C^2(\bar\Omega(0))$, $v(y,0)=0$,
$\Delta v(y,0)=0$ for $y\in \partial\Omega(0)$ and $\Delta
v(y,0)\leqslant 0$ in $\bar\Omega(0)$, then $v(y,t)\in
C^{2,1}(\bar\Omega(0)\times [0,+\infty))$ and $\Delta
v(y,t)\leqslant 0\ \textrm{for}\  y\in\Omega(0) ,\ t>0.$
\end{lem}
{\bf Proof.}  Since the initial function $v_0$ is smooth and
satisfies the consistency condition:
$$\frac{d}{\rho^2(0)}\Delta v_{0}-\frac{n\dot\rho(0)}{\rho(0)}v_0
+r v_0(1-\frac{v_0}{K})-hv_0=0\ \ \textrm{for}\ y\in
\partial\Omega(0),
$$ then the standard parabolic regularity theory \cite{LSU} shows that the
solution $v(y, t)\in C^{2, 1}(\bar\Omega(0)\times [0,+\infty))$.
Denote $w=\Delta v$, simple calculations show that it satisfies
\bess w_t \leqslant \frac{d}{\rho^2(t)}\Delta
w+\left[-\frac{n\dot\rho(t)}{\rho(t)}+r(1-\frac{2v}{K})-h\right]w.
\label{c04} \eess Taking into account the condition $\Delta
v(y,0)\leqslant 0$ we derive $w(y,0)\leqslant 0$ for $y\in
\Omega(0)$. Since $v(y,t)=0$ for $y\in
\partial\Omega(0)$, we have $$w(y,t)=\frac{\rho^2(t)}{n}\left[v_t+\frac{n\dot\rho(t)}{\rho(t)}v-rv(1-\frac{2v}{K})+hv\right](y,t)=0,\
y\in \partial\Omega(0).$$ Using the comparison principle gives that
$w(y,t)\leqslant0\ \textrm{for}\  y\in\Omega(0), \ t>0,$ which
implies that
 $\Delta v(y,t)\leqslant 0$ for $y\in\Omega(0), \ t>0.$

Let $\lambda_1$ be the principal eigenvalue of the problem
(\ref{101}) replacing $x\in\Omega$ by $y\in\Omega(0)$ then we have
the following two main theorems.

\begin{theo}\label{th31} If $0<r\leqslant \frac{d}{m^2}\lambda_1+h$,  then the solution of
problem (\ref{b12}) satisfies $v(y,t)\to 0 $ uniformly on
$\bar\Omega(0)$ as $ t\to \infty$ .
\end{theo}
{\bf Proof.} Obviously, $\hat v=0 $ is a lower solution of
(\ref{b12}). The remaining task now is to seek the upper solution of
(\ref{b12}).

To this end, define $\tilde v(y,t)$ to be the unique solution of the
problem: \bess\left\{
\begin{array}{ll}
\displaystyle \tilde v_t = \frac{d}{\rho^2(t)}\Delta\tilde
v-\frac{n\dot\rho(t)}{\rho(t)}\tilde v
 + r \tilde v(1-\frac{\tilde v}{K} )-h \tilde v ,  & y\in\Omega(0),\ t>0,
\\
\tilde v(y,t)= 0, &  y\in\partial\Omega(0), \ t>0,\\
\tilde v(y,0)=M\phi(y),  & y\in\Omega(0),
\end{array}\right.
\label{c06}\eess where $\phi $ is the corresponding eigenfunction of
$\lambda_1 $, $M$ is a positive constant. Noting the behavior of the
eigenfunction, $\phi'(y)<0$ on $\partial\Omega(0)$ for any $v_0(y)$
satisfying $v_0(y)=0$ on $
\partial\Omega$, there is $M$ such that $M\phi(y)\geqslant v_0(y)$, then $\tilde v(y,t)$ is
an upper solution of (\ref{b12}). It follows from the comparison
principle that
 \bess
  0 \leqslant v(y,t)\leqslant \tilde v(y,t),  \ y\in \Omega(0), \ t>0.
 \label{c07}
\eess Since $\Delta\tilde v(y,0)=M\Delta\phi(y)=-\lambda_1 M
\phi(y)\leqslant 0,$ it follows from Lemma 3.2 that $\Delta\tilde
v(y,t)\leqslant 0$ for $y\in \Omega(0),\ t>0$ .

On the other hand, taking into account that $\rho(t)$ tends
increasingly to $m$, $1\leqslant\rho(t)\leqslant m$ for $t\geqslant
0,$ $\tilde v(y,t)$ satisfies \bess \tilde v_t\leqslant
\frac{d}{m^2}\Delta\tilde v+r\tilde v(1-\frac{\tilde
v}{K})-h\tilde{v}. \label{c08} \eess Now consider the following
problem \bess\left\{
\begin{array}{ll}
 \displaystyle \bar v_t=\frac{d}{m^2} \Delta\bar v+r\bar v(1-\frac{\bar{v}}{K})-h\bar{v},
  & y\in \Omega(0) ,\ t>0,\\
\bar v(y,t)=0, & y\in \partial\Omega(0), \ t>0,\\
\bar v(y,0)=M\phi(y) ,  & y\in \Omega(0).
\end{array}\right.
\label{c09} \eess We may  use the comparison principle again to show
that
 $\bar v(y,t)\geqslant \tilde v(y,t)$ for $y\in \Omega(0)$ and $t>0$.
So \bess
   0\leqslant v(y,t)\leqslant \tilde v(y,t)\leqslant \bar v(y,t).
\label{c10} \eess Since that $0<r\leqslant
\frac{d}{m^2}\lambda_1+h$, we have $\bar v(y,t) \to 0$ uniformly for
$y\in \bar\Omega(0)$ as $t \to \infty$ by Theorem 1.1. Thus $v(y,t)
\to 0$ uniformly for $y\in \bar\Omega(0)$ as $ t \to \infty $.

\begin{theo} \label{th23} If $r >\frac{d}{m^2}\lambda_1+h$,  then the solution of
problem (\ref{b12}) satisfies $v(y,t)\to v^*(y)$ as \  $ t\to
\infty,$ where $v^*(y)$ is the unique positive solution of
\bes\left\{
\begin{array}{ll}
 \displaystyle-\frac{d}{m^2}\Delta v=r v(1-\frac{v}{K})-hv, & y\in \Omega(0),
\\
v(y)=0, & y\in \partial\Omega(0).
\end{array}\right.
\label{c11} \ees
\end{theo}
{\bf Proof.} Since $ \lim_{t\to\infty}\rho(t)=m$,  for any
$\varepsilon >0$, there exists a $T_0>0$, such that $m-\varepsilon
\leqslant \rho(t)\leqslant m$ for $t\geqslant T_0.$ Similarly,
 $\lim_{t\to \infty}\frac{\dot \rho(t)}{\rho(t)}=0$ implies that for
 the same $\varepsilon>0$, there exists another $T_1>0$, such that
$0\leqslant \frac{\dot\rho(t)}{\rho(t)}\leqslant \varepsilon$ for
$t\geqslant T_1.$

Set $T_*=\max \{T_0, T_1\}$ and let $\tilde v(y,t)$ denote the
solution of the following problem \bess\left\{
\begin{array}{ll}
 \displaystyle  \tilde v_t = \frac{d}{\rho^2(t)}\Delta\tilde v-\frac{n\dot\rho(t)}{\rho(t)}\tilde v
 + r\tilde v(1-\frac{\tilde{v}}{K})-h\tilde{v},  & y\in \Omega(0),\ t>T_*,
\\
\tilde v(y,t)= 0, & y\in\partial \Omega(0), \ t>T_*,\\
\tilde v(y,T_*)=M\phi(y),  & y\in \Omega(0),
\end{array}\right.
\label{c12}\eess where $M$ is a sufficiently large constant, $\phi$
is eigenfunction defined above. It follows from Lemma \ref{lem1}
that $\tilde v(y, t)$ is an upper solution of (\ref{b12}) in
$\bar\Omega(0)\times [T_*, \infty)$.

As in the proof of Theorem \ref{th31}, since that $\Delta\tilde v(y,
T_*)\leqslant 0$ in $\bar\Omega(0)$, then $\Delta\tilde v(y,
t)\leqslant 0$ in $\bar\Omega(0)\times [T_*,+\infty)$, which leads
to \bes \tilde v_t\leqslant \frac{d}{m^2}\Delta\tilde v+r \tilde
v(1- \frac{\tilde v}{K})-h\tilde{v},
 \ y\in\Omega(0),\ t>T_*.
\label{c13} \ees

Then we consider the following problem \bes\left\{
\begin{array}{ll}
 \displaystyle  v_t=\frac{d}{m^2} \Delta v+r v(1-\frac{v}{K})-hv, & y\in\Omega(0) ,\ t>T_*,
\\
v(y,t)=0, & y\in \partial \Omega(0), \ t>T_*,\\
v(y,T_*)=M\phi(y) ,  & y\in\Omega(0).
\end{array}\right.
\label{c14} \ees Clearly, the problem admits a unique solution $\bar
v(y,t)$, see \cite{P}. Moreover, since that
$r>\frac{d}{m^2}\lambda_1+h$, the result of Theorem \ref{T11} shows
that $\bar v(y,t) \to v^*(y)$ as $t\to \infty$, where $v^*(y)$ is
the unique positive solution of (\ref{c11}). Using (\ref{c13}) and
(\ref{c14}) and comparison principle yields that
\[ \tilde v(y,t) \leqslant \bar v(y,t) \ \ \mbox{for} \ \ y\in \Omega(0), \ t>T_*. \]
This implies that \bes
 \limsup_{t\to \infty} v(y,t)\leqslant v^*(y) \ \ \mbox{for}\ \ y\in \Omega(0).
 \label{c15}
 \ees

On the other hand, let $\hat v(y,t)$ be the solution of the
following problem \bess\left\{
\begin{array}{ll}
\displaystyle \hat v_t= \frac{d}{\rho^2(t)}\Delta\hat
v-\frac{n\dot\rho(t)}{\rho(t)}\hat v
 +r \hat v(1-\frac{\hat v}{K})-h\hat{v}, & y\in\Omega(0)  ,\ t>T_*,
\\
\hat v(y,t)=0, & y\in\partial\Omega(0), \ t>T_*,\\
\hat v(y,T_*)=\delta\phi(y) ,  & y\in\Omega(0),
\end{array}\right.
\label{c16} \eess where $\delta$ is a sufficiently small constant.
It is easy to see that $\hat v(y,t)$ is a lower solution of
(\ref{b12}) in $\bar\Omega(0)\times [T_*, \infty)$ if  $\delta
\phi(y) \leqslant v(y, T_*)$.

Because $\Delta\hat v(y,T_*)=-\delta\lambda_1\phi(y)\leqslant 0 $,
it follows from Lemma \ref{lem2} that $\Delta\hat v(y,t)$ $\leqslant
0$ for $y\in \bar\Omega(0),\ t\geqslant T_*$ and then
 \bess
\hat v_t\geqslant \frac{d}{(m-\varepsilon)^2}\Delta\hat v+r
\hat{v}(1-\frac{\hat v}{K})-(h+n\varepsilon)\hat v,
 \ y\in\Omega(0) ,\  t>T_*,
\label{c17} \eess since that $\rho(t)>m-\varepsilon$ and
$\frac{\dot\rho(t)}{\rho(t)}\leqslant\varepsilon$ for $t\geqslant
T_*$.

Now consider the following problem
 \bes\left\{
\begin{array}{lll}
\displaystyle \hat v_t= \frac{d}{(m-\varepsilon)^2}\Delta\hat v
 +r\hat v(1- \frac{\hat{v}}{K})-(h+n\varepsilon)\hat v,
 & y\in\Omega(0) ,\  t>T_*,
\\
\hat  v(y,t)=0, & y\in\partial\Omega(0), \ t>T_*,\\
\hat v(y,T_*)=\delta\phi(y),  & y\in\Omega(0).
\end{array}\right.
\label{c18} \ees Similarly (\ref{c18}) admits a unique positive
solution, denoted by $\hat v_\varepsilon(y,t)$. Using comparison
principle yields that $\hat v_\varepsilon(y,t)\leqslant \hat
v(y,t)$. Since $r>\frac{d}{m^2}\lambda_1+h$, we can choose
$\varepsilon>0$ sufficiently small such that
$r>\frac{d}{(m-\varepsilon)^2}\lambda_1+h+n\varepsilon$. we then
have $\hat v_\varepsilon(y,t) \to \hat v_\varepsilon^*(y)$ as $t \to
\infty$, where $\hat v_\varepsilon^*(y)$ is the unique positive
solution of \bess\left\{
\begin{array}{ll}
\displaystyle -\frac{d}{(m-\varepsilon)^2}\Delta\hat v
 =(r-n\varepsilon-h)\hat v-\frac{r}{K}\hat v^2,  & y\in\Omega(0),
 \\
\hat v(y)=0, & y\in\partial\Omega(0).
\end{array}\right.
\label{c19} \eess Therefore we have \bes
 \liminf_{t\to \infty} v(y,t)\geqslant \hat v_\varepsilon^*(y)\ \ \textrm{for} \  y\in \Omega(0).
   \label{c20}
  \ees
By the continuous dependence of $\hat v_\varepsilon^*(y)$ on
$\varepsilon$, we can easily see that $\hat v_\varepsilon^*(y) \to
v^*(y)$ as $\varepsilon \to 0^+$. Thus, (\ref{c15}), (\ref{c20}) and
the arbitrariness of $\varepsilon$ imply $v(y,t)\to v^*(y)$
uniformly on $\bar\Omega(0)$ as $ t\to \infty$, where $v^*(y)$
satisfies problem (\ref{c11}). Thus Theorem \ref{th23} is proved.

\section{Numerical results}\setcounter{equation}{0}
As we know, the presented harvest model is of practical interest
only for one or two space variables ($n=1,2$). In this section, we
restrict $\Omega$ to one dimensional space. Numerical simulation is
carried out to illustrate the results obtained in previous sections
and to explore the long-time behavior of solutions on growing
domain. The essential idea of the numerical calculation is to
transform the growing domain to a fixed domain as done in equations
(\ref{b12}). The consequence is that the diffusivity changes from
being time-independent on the growing domain to being time-dependent
on the fixed domain. On unvaried domain, there are a handful of
well-developed numerical methods that can be used directly.

Firstly, regarding the domain growth, we choose
$\Omega(t)=[0,x(t))=[0,\rho(t)y)$, where the parameters $k=1, m=2$.
It is easy to see that domain grows from initial size $\rho(0)=1$ to
the final size $\rho(\infty)=m=2$. We take initial function as
$u(y,0)= \sin(y),\ y\in(0,1)$ corresponding to $t=0$.

Then, we take the values of the parameters in system (\ref{b12}) as
follows: \bes d=0.9,\ r=2,\ K=4,\ h=0.5.
 \label{d01}
 \ees
In the case of the interval $(0,1)$, $\lambda_1=\pi^2$. Under the
set of parameters in (\ref{d01}), we have
$0<r<d(\frac{\pi}{m})^2+h\approx2.72$. By Theorem \ref{th31} we know
that the solution of problem (\ref{b12}) satisfies $v(y,t)\to 0 $
uniformly on $[0,1]$ as $ t\to \infty$. Then we have $u(x,t)\to 0$
uniformly on any compact subset of $[0,2]$ as $ t\to \infty$. This
is shown in Figure 1-(a), where the process of domain growth is
presented in the left figure.
\begin{figure}[hbtp]\label{f1}
 \begin{center}
 \includegraphics[width=0.9\textwidth]{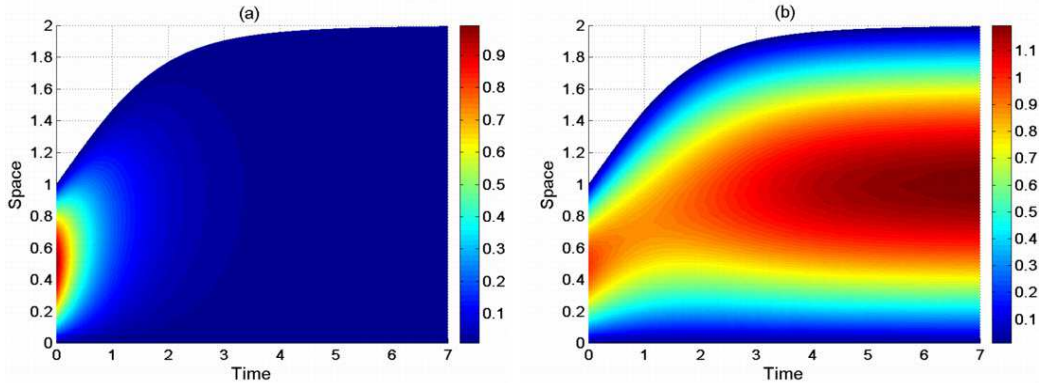}
  \caption{Asymptotic behavior of the solution to system (\ref{b12}). (a): $r<d(\frac{\pi}{m})^2+h$, the parameters are given in
  (\ref{d01}).
  (b): $r>d(\frac{\pi}{m})^2+h$, the parameters are given in (\ref{d02}).}
\end{center}
 \end{figure}
In order to illustrate the result in Theorem \ref{th23}, we choose
an other set of the parameters
 \bes
 d=0.9, \ r=4,\ K=4,\ h=0.5 .
 \label{d02}\ees
In this case, $r>d(\frac{\pi}{m})^2+h\approx 2.72$, we know that
solution $v(y,t)$ of problem (\ref{b12}) asymptotically converges to
the steady state $v^*(y)$ by Theorem \ref{th23}. Figure 1-(b) shows
that on the growing domain $(0,x(t))$, the solution
$u(x,t)(=u(\rho(t)y,t)=v(y,t))$ asymptotically converges to the
steady state $v^*(y)$. From the numerical simulation, we can see
that most of the individuals aggregate around the center of domain
as time increases.

On the other hand, from the Theorem \ref{th23} we know if
$r>d(\frac{\pi}{m})^2+h$,  the species will tend to extinction when
the harvesting rate $h$ increases. This can be seen from Figure 2,
where we take $h=1$ and $1.5$, respectively. The other parameters
are same as in (\ref{d02}).
\begin{figure}[hbtp]\label{f2}
 \begin{center}  \includegraphics[width=0.9\textwidth]{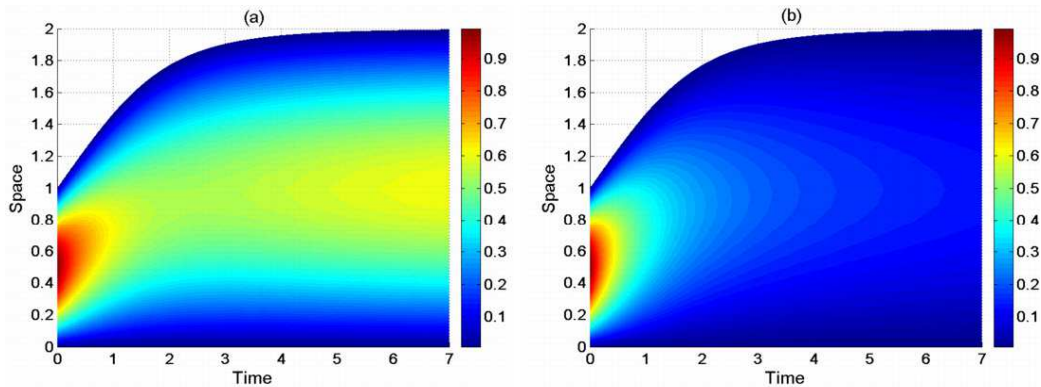}
  \caption{Asymptotic behavior of the solution to system (\ref{b12}) when $r>d(\frac{\pi}{m})^2+h$. (a)$h=1$, (b)$h=1.5$, the other parameters are same as
  in (\ref{d02}). }
\end{center}
 \end{figure}

\section{Conclusion}
Recently, domain growth is an interesting topic which has attracted
a lot of attention. However, most existing results on the long time
behaviors of the solutions were investigated through numerical
simulations. In this paper, we succeeded achieving the global
stability of the solution to a harvest single species logistic model
with an isotropic domain growth was studied under Dirichlet boundary
condition via upper and lower solutions. We first developed model
and then verified the comparison principle which is fundamentally
important in studying the asymptotical behavior of temporal
solutions to problem (\ref{b12}). Then asymptotic behavior of
solutions was investigated by approach of upper and lower solutions.
Our results show that $v(y,t)(=u(\rho(t)y, t))$  converges to $0$ if
$r\leqslant \frac{d}{m^2}\lambda_1+h$ or to the nonnegative steady
state solution if $r> \frac{d}{m^2}\lambda_1+h$. Finally we show
that numerical simulations are consistent with our analytical
results. Of course, this method allows to obtain the asymptotic
estimates for the more general growth functions, which are monotone
and continuous differentiable on $[0,+\infty)$. 

Ecologically speaking, the results imply that the growth of domain has a positive effect on the asymptotic stability of positive steady state solution and a negative effect on the asymptotic stability of the trivial solution. Conversely, the impact of the harvesting rate is opposite: a negative effect on the asymptotic stability of positive steady state solution and a positive effect on the asymptotic stability of the trivial solution. In other words, if the harvesting rate is large, the species is more at risk for extinction.

\bibliographystyle{my}
\small\bibliography{reference}

\end{document}